\documentclass[10pt,a4paper]{amsart}
\usepackage{amsmath}
\usepackage{amssymb}
\usepackage{amscd}
\input xypic

\newtheorem{Thm}{Theorem}[section]
\newtheorem{Lem}[Thm]{Lemma}
\newtheorem{Cor}[Thm]{Corollary}
\newtheorem{Prop}[Thm]{Proposition}

\newtheorem{Rmk}[Thm]{Remark}

\def\2m#1#2#3#4{\left(\matrix{#1&#2\cr#3&#4\cr}\right)}
\def\vect#1#2{\left(\begin{matrix}#1\cr #2\cr\end{matrix}\right)}

\def\C{\mathcal C}
\def\D{\mathcal D}

\def\F{\mathcal F}
\def\G{\mathcal G}
\def\H{\mathcal H}
\def\I{\mathcal I}
\def\J{\mathcal J}
\def\K{\mathcal K}

\def\M{\mathcal M}

\def\P{\mathcal P}

\def\T{\mathcal T}

\def\Y{\mathcal Y}

\def\2matr#1#2#3#4{\left(
      \begin{array}{cc}
#1&#2\\
#3&#4\\
      \end{array}
      \right)}
\def\matrQ#1#2#3#4{\left(
      \begin{array}{cc}
#1&#2\\
#3&#4\\
      \end{array}
      \right)}

\def\matrHH#1#2#3{\left(
      \begin{array}{ccc}
#1&#2&#3\\
      \end{array}
      \right)}
\def\matrQQ#1#2#3#4#5#6#7#8#9{\left(
      \begin{array}{ccc}
#1&#2&#3\\
#4&#5&#6\\
#7&#8&#9\\
      \end{array}
      \right)}
\def\limprojder#1#2#3
{{\mathop{\varprojlim}_{#2}}^{(#1)}#3}

\def\limproj#1#2
{\mathop{\varprojlim}
_{#1}#2}

\def\liminjder#1#2#3
{\mathop{{\varinjlim}^{(#1)}}
_{#2}#3}

\def\liminj#1#2
{\mathop{\varinjlim}
_{#1}#2}

\begin{document}

\author{Arne B. Sletsj\o{}e}

\address{
Matematisk institutt,\newline
University of Oslo,\newline
Pb. 1053, Blindern,\newline
N-0316 Oslo, Norway}

\email{arnebs@math.uio.no}


\begin{abstract}
    We study finite dimensional representations of the projective modular
group. Various explicit dimension formulas are given. 
\end{abstract}

\title{Finite dimensional representations of the projective 
modular group}
\maketitle


\section{Introduction}\label{Introduction}
\subsection{Finite-dimensional modules}\label{Finite-dimensional modules}
Let $k$ be an algebraically closed field of characteristic zero. 
An $n$-dimensional left module of a $k$-algebra
$A$ is a vector space $M$ of dimension $n$ over $k$, equipped with a structure map
 $\rho\in Hom_{k\textrm{-}alg}(A,End_{k}(M))$.
The (left) module $M$ is simple if it has no non-trivial
submodules, or equivalently, if the structure map $\rho$ is surjective.
A module which does not split into a direct sum of submodules is
called indecomposable. To each indecomposable
module $M$ we associate a unique, up to permutations of the components,
semi-simple module $\overline{M}$, given as the direct sum of the composition
factors of $M$.
\par
Two structure maps $\rho, \rho^{\prime}$ define equivalent modules if they are
conjugates, i.e. there exists an invertible $n\times n$-matrix $P$ such
that for every $a\in A$, $\rho^{\prime}(a)=P^{-1}\rho (a)P$. The
set of conjugation classes of modules has in general no structure of an algebraic
variety, but restricting to the semi-simple modules of dimension $n$ there exist
an algebraic
quotient, denoted by $\M_{n}$. The variety $\M_{n}$ can also be considered as 
a quotient of the variety of all $n$-dimensional $A$-modules under the weaker
equivalence relation, where
two modules are equivalent if they have conjugate composition factors. The
Zariski open subset of $\M_{n}$ of simple modules is denoted $Simp_{n}(A)$.
This subvariety is independent of the choice of equivalence relation. 
\subsection{The modular group}\label{The modular group}

The modular group $SL(2,\mathbf{Z})$ is generated by the two matrices 
\begin{displaymath}
   u=\matrQ11{-1}0,\qquad  v=\matrQ01{-1}0
\end{displaymath}
where $\overline{u}^3=\overline{v}^2=I$ in the projectivised group 
$PSL(2,\mathbf{Z})\simeq SL(2,\mathbf{Z})/\{\pm I\}$. 
There are no relations 
between the two generators and $PSL(2,\mathbf{Z})$ is clearly
isomorphic to the free product $\mathbf{Z}_{3}\ast\mathbf{Z}_{2}$. 
\par
Finite-dimensional representations of a group $G$ correspond 
bijectively to finite-dimensional modules of the group 
algebra $k[G]$. Thus the representation theory of $PSL(2,\mathbf{Z})$ 
is equivalent to the representation theory of the group algebra
\begin{displaymath}
    A=k[\mathbf{Z}_{3}\ast\mathbf{Z}_{2}]\simeq
    k[\mathbf{Z}_{3}]\ast k[\mathbf{Z}_{2}] \simeq
     k\langle x,y\rangle/(x^3-1,y^2-1)
\end{displaymath}
The algebra $k[\mathbf{Z}_{p}]$ is hereditary, and 
free products of hereditary algebras are hereditary. Thus the 
algebra $A$ is hereditary.
\par
Let $M$ be some finite dimensional representation of the non-commutative model
$B=k\langle x,y\rangle/(x^3-y^2)$ of the ordinary cusp. The inclusion of the
center
$Z(B)\simeq k[x^3=y^2]\hookrightarrow B$ induces a fibration
\begin{displaymath}
Simp_{n}(B)\longrightarrow \mathbf{A}^1
\end{displaymath}
for which the fibres outide the origin are constant, isomorphic to
$Simp_{n}(A)$.
\subsection{Low-dimensional modules}\label{Low-dimensional
modules}

We shall give an explicit description of the low-dimensional modules of 
$A=k\langle x,y\rangle/(x^3-1,y^2-1)$. For further details we refer to
\cite{TW}.
\par
The variety 
$Simp_{1}(A)$ consists of 6 points, denoted $k(\omega,\alpha)$, and defined by
$x\mapsto 
\omega$, $\omega^3=1$ and $y\mapsto \alpha=\pm 1$. 
We have 
$Ext_{A}^1(k(\omega,\alpha),k(\omega^{\prime},\alpha^{\prime}))\simeq 
k$ if $\omega\neq \omega^{\prime}$ and $\alpha\neq \alpha^{\prime}$,  
and zero elsewhere.
\par
The variety of two-dimensional simple modules 
 $Simp_{2}(A)$
has three disjoint components, 
and each component is an affine line with two distinct closed points
removed. 
We observ that there exist non-conjugate indecomposable modules 
 $E_{1}$ and $E_{2}$, satisfying $Ext^1_{A}(E_{1},E_{2})\neq 0$.
Thus it is not possible to complete $Simp_{2}(A)$ to an ordinary commutative
scheme, under the strong conjugation relation.
To illustrate, for $s\in k$ we have a 
2-dimensional $A$-module
\begin{displaymath}
    M_{s}:\quad x\mapsto\matrQ110{\omega}, \qquad y\mapsto\matrQ10s{-1}
\end{displaymath}
where $\omega^3=1$, $\omega\neq 1$. For $s\neq 0, 2(\omega-1)$ the
modules are simple, while $M_{0}$ and $M_{2(\omega-1)}$ are indecomposable, 
non-simple, corresponding to extensions of $k(\omega,-1)$ by 
$k(1,1)$, respectively $k(\omega,1)$ by $k(1,-1)$. The trivial extensions are
the respective associated semi-simple modules. In fact all simple modules have
a representation of this form.
\par
There is yet another family, given by
\begin{displaymath}
    N_{t}:\quad x\mapsto\matrQ101{\omega}, \qquad y\mapsto\matrQ1t0{-1}
\end{displaymath}
For $t\neq 0, 2(\omega-1)$ the
modules are simple and $N_{t}\simeq M_{t}$. As above the two exeptional modules
$N_{0}$ and $N_{2(\omega-1)}$ are indecomposable and 
non-simple, but not isomorphic to any $M_{s}$. The module $N_{0}$ corresponds to
an extension of $k(1,-1)$ by 
$k(\omega,1)$, the opposite direction of $M_{0}$, and $Ext^1_{A}(M_{0},N_{0})\neq
0$.
\par
For $n=3$ the situation is as follows. Two matrices $P$ and 
$Q$ are given by
\begin{displaymath}
    P=\matrQQ{\lambda_{1}}
    {\frac{\lambda_{1}\lambda_{3}}{\lambda_{2}}+\lambda_{2}}
    {\lambda_{2}}0{\lambda_{2}}{\lambda_{2}}00{\lambda_{3}}
    \qquad
    Q=\matrQQ{\lambda_{3}}00{-\lambda_{2}}{\lambda_{2}}0
    {\lambda_{2}}
    {-\frac{\lambda_{1}\lambda_{3}}{\lambda_{2}}-\lambda_{2}}
    {\lambda_{1}}
\end{displaymath}
where $(\lambda_{1}\lambda_{2}\lambda_{3})^2= 1$. The variety $Simp_{3}(A)$ has
two components, given by $\lambda_{1}\lambda_{2}\lambda_{3}= \pm 1$. For the
component $\lambda_{1}\lambda_{2}\lambda_{3}= 1$ we put
\begin{displaymath}
    x\mapsto PQ=\matrQQ00{\lambda_{1}\lambda_{2}}
    0{-\lambda_{1}\lambda_{3}}{\lambda_{1}\lambda_{2}}
    {\lambda_{2}\lambda_{3}}{-\frac{\lambda_{1}\lambda_{3}^2}
    {\lambda_{2}}-\lambda_{2}\lambda_{3}}
    {\lambda_{1}\lambda_{3}}
    \qquad
    y\mapsto PQP=\matrQQ0010{-1}0100
\end{displaymath}
defining a 3-dimensional module of $A$.
A consequence of the Cayley-Hamilton Theorem is that since $x^3=y^2=1$, the
images of the 9 monomials 
\begin{displaymath}
    1,x,y,x^2,xy,yx,x^2y,xyx,x^2yx
\end{displaymath}
form a linear basis for $M_{3}(k)$ if and only if the module is simple.
A computation shows that this set of matrices is linearily independent if and
only if
\begin{displaymath}
    (\lambda_{1}^3+1)(\lambda_{2}^3+1)(\lambda_{3}^3+1)\neq 0
\end{displaymath}
One can show that in this case the trace ring (\cite{JLS}) is generated by
\begin{align*}
    t_{xy}&=\frac{1}{\lambda_{1}}+\frac{1}{\lambda_{2}}+\frac{1}{\lambda_{3}}\\
    t_{xyx}&=\lambda_{1}+\lambda_{2}+\lambda_{3}\\
\end{align*}
For any choice of cube root $\lambda$ of -1 we get the relation
\begin{displaymath}
    t_{xy}-\lambda t_{xyx}+2\lambda^2=0
\end{displaymath}
corresponding to a straight line in the $(t_{xy}, t_{xyx})$-plane. 
Thus the non-simple locus is given by the three lines
\begin{displaymath}
    t_{xy}-\lambda t_{xyx}+2\lambda^2=0
\end{displaymath}
for $\lambda^3=-1$. The three intersections of the three lines correspond to
semi-simple modules which are sums of three one-dimensional simple modules.
\par
Similar analysis for $Simp_{n}(A)$ for $n=4,5$, can be found in \cite{TW}.
\par

\subsection{Main results}\label{Main results}
In this paper we investigate the finite dimensional modules
of the projective modular group $PSL(2,\mathbf{Z})$. Unfortunately the litterature
is somewhat
limited on this subject.
We have mentioned the paper of Tuba and Wenzl (\cite{TW}).
Other important papers are due to Westbury (\cite{W}) and Le Bruyn and
Adriaenssen (\cite{A-LB}). These papers give nice desriptions
of the finite-dimensional modules, especially in low dimensions, i.e. up to
dimension 5.
Le Bruyn and Adriaenssens also give modules of dimension $6n$. Our aim in this
paper is to give more explicit results on the varieties of modules
of arbitrary dimensions.
\par
The paper is organised as follows. 
In section \ref{Quiver theory} we recall a result of Wesbury \cite{W}, proving
a dimension formulae for the variety of simple $n$-dimensional $A$-modules of
given dimension vector
$\alpha=(\alpha_{1},\alpha_{2},\alpha_{3};
\alpha_{1}^{\prime},\alpha_{2}^{\prime})$, denoted $Simp_{\alpha}(A)$.
Suppose $\alpha$ is a Schur root, i.e. there exist a module $M$ of dimension
vector $\alpha$ such that $End_{A}(M)\simeq k$.
Westburys result tells us that in that case 
\begin{displaymath}
dim\, Simp_{\alpha}(A)=n^2-\sum_{i=1}^3 \alpha_{i}^2
-\sum_{j=1}^2 (\alpha^{\prime}_{j})^2 +1
\end{displaymath}
In section \ref{Deformation theory} we look at the locus $\D_{\alpha}$ of
weak equivalence classes of non-simple modules, represented by the associated
semi-simple module. Using deformation theory we prove that 
$\overline{E}\in D_{\alpha}$, corresponding to an extension of two non-isomorphic
simple modules of dimension vector $\beta$ and $\gamma$, is a smooth point of 
a component of $D_{\alpha}$ of codimension $-\langle\beta, \gamma, \rangle$.
\par
In section \ref{Maximally iterated extension} we show that there exist 
$\frac{1}{2}(n+d+1)$ distinct
points of $\D_{\alpha}$ corresponding to the maximally iterated extensions in the
given
component. Each of these points in the weak equivalence relation corresponds
to a set of strong equivalence classes. In section \ref{Parametrisation}
we study this collection more carefully. The weak equivalence classes split
up in several equivalence classes under the strong relation. The number of
strong classes depends not only on the dimension vector $\alpha$ and the
composition factors, but also on the order in which the iterated extensions
are composed, described by the representation graph $\Gamma$. We show that the
non-commutative scheme $Ind_{\Gamma}(A)$ can be described as
a quotient of a certain vector space by a group of conjugations.
\par
In section \ref{Modular Forms} we 
consider the component of $\M_{n}$ of highest dimension, and inside this the
subvariety of non-simple modules. It turns out that the minimal codimension of
this non-simple locus is the same as the dimension of the space of modular
forms of weight $2n$.
\par

\section{The Quiver $K(3,2)$}\label{Quiver theory}
\subsection{Quiver representation}\label{Quiver representations}

Let $K(3,2)$ be the quiver with 
vertices $v_{1},v_{2},v_{3}$ and $w_{1},w_{2}$ and with an 
arrow from $v_{i}$ to $w_{j}$ for $i=1,2,3$ and $j=1,2$. The variety of
equivalence classes of (semi-)simple
representations of $\mathbf{Z}_{3}\ast\mathbf{Z}_{2}$ corresponds to the
variety of equivalence classes of $\theta$-(semi)stable representations of
$K(3,2)$ (\cite{W}),
where $\theta$ is the character assigning to the dimension vector
$\alpha=(\alpha_{1}, \alpha_{2}, \alpha_{3},\alpha^{\prime}_{1},
\alpha^{\prime}_{2})$ the dimension difference 
\begin{displaymath}
\theta(\alpha)=(\alpha^{\prime}_{1}+\alpha^{\prime}_{2})-
(\alpha_{1}+\alpha_{2}+\alpha_{3})
\end{displaymath}
\par
The correspondance goes via the category
$\C(\mathbf{Z}_{3}\ast\mathbf{Z}_{2})$, whose objects are triples
$(M,N,\phi)$, where $M$ is a left $k[\mathbf{Z}_{3}]$-module, $N$ is a left
$k[\mathbf{Z}_{2}]$-module and $\phi:N\rightarrow M$ is a $k$-linear map. To a
semi-simple representation $M$ of $\mathbf{Z}_{3}\ast\mathbf{Z}_{2}$ we associate the triple
$(M,M,id)\in \C(\mathbf{Z}_{3}\ast\mathbf{Z}_{2})$, 
where the $\mathbf{Z}_{3}$- and $\mathbf{Z}_{2}$-module
structures of $M$ are induced by the two inclusions 
\begin{displaymath}
\mathbf{Z}_{3}\hookrightarrow \mathbf{Z}_{3}\ast
\mathbf{Z}_{2}\hookleftarrow \mathbf{Z}_{2}
\end{displaymath}
The triple $(M,M,id)$ is mapped to the $K(3,2)$-module $M\oplus M$ where
$v_{i}(m,m^{\prime})=(e_{i}m,0)$ for an idempotent
$e_{i}=\frac{1}{3}(1+\omega^{i}
x+\omega^{2i} x^2)\in k[\mathbf{Z}_{3}]\simeq k[x]/(x^3-1)$, 
$w_{j}(m,m^{\prime})=(0,f_{j}m^{\prime})$, with
$f_{j}=\frac{1}{2}(1+(-1)^jy)\in k[\mathbf{Z}_{2}]\simeq k[y]/(y^2-1)$ and 
$\phi_{ij}(m,m^{\prime})=(0,f_{j}e_{i}m)$.
Here $i=1,2,3$, $j=1,2$ and $\omega$ is a primitive cube root of unity.
Thus the eigenspaces of the generators $x$ and $y$ correspond to the components
of the module $M\oplus M$ over the quiver.
\par

Let $M$ be an indecomposable $A$-module with dimension vector denoted by
$\alpha=(\alpha_{1},\alpha_{2},\alpha_{3};\alpha^{\prime}_{1},
\alpha^{\prime}_{2})$ such that 
$\vert \alpha\vert =\alpha_{1}+\alpha_{2}+\alpha_{3}
=\alpha^{\prime}_{1}+\alpha^{\prime}_{2}=n$. 
The dimension 
of  $Ext_{A}^1(M,M)$ is given by the formula (\cite{Sch})
\begin{displaymath}
    dim\,Ext_{A}^1(M,M)=dim \,Hom_{A}(M,M)-\langle \alpha,\alpha\rangle
\end{displaymath}
where 
\begin{displaymath}
    \langle \alpha, \beta\rangle=\sum_{i=1}^3\alpha_{i}\beta_{i}
    +\sum_{j=1}^2\alpha^{\prime}_{j}\beta^{\prime}_{j}
    -(\sum_{i=1}^3\alpha_{i})(\sum_{i=1}^3\beta_{i})
\end{displaymath}

Thus if $End_{A}(M)\simeq k$, we get
\begin{displaymath}
    d=dim\,Ext_{A}^1(M,M)=1+n^2-\sum_{i=1}^3 {\alpha_{i}}^2-\sum_{j=1}^2 
    (\alpha^{\prime}_{j})^2
\end{displaymath}

\subsection{Westburys theorem.}\label{Westburys theorem}

The following theorem is due to Westbury (\cite{W}).

\begin{Thm}\label{Westbury}
    Let $\alpha=(\alpha_{1},\alpha_{2},\alpha_{3};
    \alpha^{\prime}_{1},\alpha^{\prime}_{2})$ be a dimension vector of 
    non-negative integers, satisfying
    \begin{itemize}
\item[i)] $\alpha_{1}+\alpha_{2}+\alpha_{3}
=\alpha^{\prime}_{1}+\alpha^{\prime}_{2}=n$
\item[ii)] $\alpha_{i}+\alpha^{\prime}_{j}\le n$ for all $i=1,2,3$ and $j=1,2$
    \end{itemize}
Then the component $Simp_{\alpha}(A)$ of
$Simp_{n}(A)$ corresponding to the dimension 
    vector $\alpha$ is a non-empty smooth affine variety of
dimension 
    \begin{displaymath}
d=dim\, Simp_{\alpha}(A)=1+n^2-\sum_{i=1}^3 {\alpha_{i}}^2
-\sum_{j=1}^2 (\alpha^{\prime}_{j})^2 
    \end{displaymath}
\end{Thm}

\begin{Cor}
    The dimension of $Simp_{n}(A)$, i.e. the largest dimension of its 
    components is given by
    \begin{displaymath}
dim\, Simp_{n}(A)=
\begin{cases}
6m^2+2sm+s-1& n=6m+s, s=1,\ldots,5, m\ge 0\\
6m^2+1& n=6m, m\ge 1\\
\end{cases}
    \end{displaymath}
\end{Cor}

\begin{proof}
    It is easy to see that maximal dimension is obtained when the 
    partition of $n$ is most evenly distributed. Using this general 
    principle the computation of the dimensions is straight-forward.
\end{proof}

\begin{Rmk}
    The generating function of the sequence of maximal dimensions is given by
    \begin{displaymath}
\sum_{n=1}^{\infty} dim\, (Simp_{n}(A))t^{n}
=\frac{t^2+2t^6-2t^7+t^8}{(1-t)^2(1-t^6)}
    \end{displaymath}
\end{Rmk}

\section{Deformation theory}\label{Deformation theory}

\subsection{Extensions for the modular group}\label{Extensions for the modular
group}
The variety $\M_{\alpha}$ of finite dimensional modules
of $A$ of given dimension vector $\alpha$, up to weak equivalence, where two
modules are equivalent if they have conjugate composition factors, is a coarse
moduli space.
The $\alpha$-component of the variety of simple modules $Simp_{\alpha}(A)\subset
\M_{\alpha}$ is a Zariski open subset of $\M_{\alpha}$. 
The complement in $\M_{\alpha}$
of weak equivalence classes of non-simple modules
is denoted $\D_{\alpha}$.
The dimension of the tangent space of $\M_{\alpha}$ at an indecomposable
module $M$ is given in section
\ref{Quiver theory} by the formula for
the dimension of $Ext_{A}^1(M,M)$. The ring $A$ is hereditary so there are no
obstructions and the coarse moduli is a smooth space of the given dimension.
\par
Let $E$ be an indecomposable, non-simple representation and suppose 
\begin{displaymath}
    0\rightarrow W\longrightarrow E\longrightarrow 
    V\rightarrow 0
\end{displaymath}
is a non-split exact sequence, presenting the indecomposable module $E$ as an
extension of simple $A$-modules $V$ and $W$.
At the point $\overline{E}\in\D_{\alpha}$ the tangent space of
$\M_{\alpha}$ contains directions corresponding to simple deformations as well
as directions correponding to non-simple modules. To study the non-simple
deformation we consider deformations of $E$ under the extension constrain,
i.e. liftings of $E$ that preserves the exact sequence presentation. This is
the topic of the next subsections.

\subsection{Deformations of presheaves}\label{Deformations of presheaves}

Let $\underline{a}$ be the category of commutative pointed artinian 
$k$-algebras, i.e. objects are commutative artinian $k$-algebras $R$ with 
a fixed quotient $\pi:R\rightarrow k$ and morphisms preserving the point. 
\par
Let $R\in obj(\underline{a})$ and let $\F:\Lambda\rightarrow 
A\textrm{-}\textbf{mod}$ be a presheaf of left $A$-modules on a poset
$\Lambda$, i.e. $\F$ is a contravariant functor on $\Lambda$ with values in
the category of left $A$-modules. 
A \textbf{lifting} of the presheaf 
$\F$ to $R$ is a presheaf $\F_{R}$ of left $A\otimes_{k}R^{op}$-modules
\begin{displaymath}
    \F_{R}:\Lambda\longrightarrow 
    A\otimes_{k}R^{op}\textrm{-}\textbf{mod}
\end{displaymath}
together with an isomorphism of presheaves
\begin{displaymath}
    \eta:\F_{R}\otimes_{R}k\buildrel\cong\over\longrightarrow\F
\end{displaymath}
and such that $\F_{R}\simeq\F\otimes_{k}R$ as presheaves of right $R$-modules.

Two liftings $\F_{R}^{\prime}$ and $\F_{R}^{\prime\prime}$ 
are said to be \textbf{equivalent} if there exists an isomorphism 
\begin{displaymath}
    \phi:\F_{R}^{\prime}\longrightarrow \F_{R}^{\prime\prime}
\end{displaymath}
of presheaves such that $\eta^{\prime\prime}\circ(\phi\otimes 
id_{k})=\phi\circ \eta^{\prime}$.
Denote by 
\begin{displaymath}
    Def_{\F}(R)=\{liftings\,\,of\,\,\F\,\,to\,\,R\}/\backsim
\end{displaymath}
the set of equivalence classes of liftings of $\F$ to $R$. We shall 
refer to such equivalence classes as \textbf{deformations} of $\F$ to $R$. This
construction is functorial, defining a covariant functor
\begin{displaymath}
    Def_{\F}:\underline{a}\longrightarrow \underline{sets}
\end{displaymath}
\par
One can show that there exist an element 
$\psi=(\psi^2,\psi^1,\psi^0)$
representing a class of a certain cohomology group
$Ext^2_{\Lambda}(\F,\F\otimes ker\,\pi)$
such that 
$[\psi]=0$ in $Ext^2_{\Lambda}(\F,\F\otimes ker\,\pi)$ is equivalent to the
existence of a lifting of $\F$ to $R$.
If $[\psi]=0$, then the set of deformations is given as a principal homogenous
space over $Ext^1_{\Lambda}(\F,\F\otimes ker \,\pi)$. 
\par
Applying the deformation functor to the $k$-algebra $R=k[x]/(x^2)$, we have
$ker\,\pi\simeq k$, and the tangent space of the deformation functor
$Def_{\F}$ is given by $Ext^1_{\Lambda}(\F,\F)$. 
In the next section we shall study $Ext^{i}_{\Lambda}(\F,\G)$ for $i\ge 0$ in
more details.

\subsection{A double complex}\label{A double complex}

Let $\F,\G:\Lambda\rightarrow A\textrm{-}\textbf{mod}$ be presheaves of 
left $A$-modules on a poset
$\Lambda$.  Define a covariant functor 
\begin{displaymath}
Hom_{k}(\F,\G):Mor(\Lambda)\longrightarrow A\textrm{-}bimod
\end{displaymath}
by $Hom_{k}(\F,\G)(\lambda^{\prime}\le\lambda)=
Hom_{k}(\F(\lambda),\G(\lambda^{\prime}))$.
There is a double complex
\begin{displaymath}
    K^{\bullet\bullet}=D^{\bullet}(\Lambda,C^{\bullet})
\end{displaymath}
given by
\begin{displaymath}
    K^{p,q}=\prod_{\Lambda_{0}\le..\le\lambda_{p}}
    C^q(A,Hom_{k}(\F(\lambda_{p}),\G(\lambda_{0})))
\end{displaymath}
where $C^q(A,\H)$ is the Hochschild cochain complex, 
with Hochschild differential
\begin{displaymath}
    d:K^{p,q}\longrightarrow K^{p,q+1}
\end{displaymath}
and where
\begin{displaymath}
    \delta:K^{p,q}\longrightarrow K^{p+1,q}
\end{displaymath}
is the differential of the $D^{\bullet}$-complex given by
\begin{align*}
    \delta\psi(\lambda_{0}\le\ldots\le\lambda_{{p+1}})=
    &\G(\lambda_{0}\le\lambda_{1})
    \psi(\lambda_{1}\le\ldots\le\lambda_{p+1})\\
    &+\sum_{i=1}^p(-1)^i
    \psi(\lambda_{0}\le\ldots\le\hat{\lambda_{i}}\le\ldots\le\lambda_{p+1})\\
    &+(-1)^{p+1}
    \psi(\lambda_{0}\le\ldots\le\lambda_{p})
    \F(\lambda_{p}\le\lambda_{p+1})\\
\end{align*}

Let $Tot(K^{\bullet\bullet})$ be the total complex 
with total differential  $\partial=d+(-1)^p\delta$, 
satisfying $\partial^2=0$. Define 
\begin{displaymath}
    Ext^n_{\Lambda}(\F,\G)=H^n(Tot(K^{\bullet\bullet}))\qquad n\ge 0 
\end{displaymath}
As indicated in the last section we are interested in $Ext^n_{\Lambda}(\F,\F)$
for $n=0,1,2$. Our main tool will be the spectral sequence of the next
proposition. Denote by $s,t$ the source and target functors on $Mor(\Lambda)$.
\begin{Prop}\label{SQ1}
    Let $A$ be a $k$-algebra, $\Lambda$ a poset and $\F, \G$ two presheaves 
    of left $A$-modules on $\Lambda$. Then there exists a first quadrant 
    spectral sequence 
    \begin{displaymath}
           E_{2}^{p,q}=\limprojder{p}{Mor(\Lambda)}{Ext^q_{A}(\F(t),\G(s))}
    \end{displaymath}
    converging to the total cohomology $Ext_{\Lambda}^{p+q}(\F,\G)$.
\end{Prop}

\begin{proof} 
    The spectral sequence is obtained by computing cohomology of the 
    double complex $K^{\bullet\bullet}$ with respect to the 
    Hochschild differential and then with respect to the 
    $D^{\bullet}$-differential.
\end{proof}

\subsection{Application to the poset $\Lambda_{2}=\{0\le 1\le 2\}$}
A short excact sequence 
    \begin{displaymath}
\F:\qquad 0\rightarrow W\buildrel{\epsilon}\over{\longrightarrow}
E\buildrel{\pi}\over{\longrightarrow} V\rightarrow 0
    \end{displaymath}
of left $A$-modules 
can be considered as a presheaf on $\Lambda_{2}$ with certain obvious properties.
We put $\F(2)=W$, $\F(1)=E$, $\F(0)=V$, $\F(0\le 1)=\pi$ and $\F(1\le
2)=\epsilon$. The smooth moduli of deformations of the $A$-module $E$ has
a tangent space $\T_{\M,E}$ given by $Ext_{A}^{1}(E,E)$ and the subspace
$\T_{\M,E}^{\F}\subset \T_{\M,E}$ coresponding to directions in $\T_{\M,E}$
subject to the extension constrain, is given by the image of the projection map 
\begin{displaymath}
Ext_{\Lambda}^{1}(\F,\F)\longrightarrow Ext_{A}^1(E,E)
\end{displaymath}
Denote this image by $Ext_{A}^1(E,E)_{0}$. 
Observe that the longest sequence of non-trivial elements of 
$Mor(\Lambda_{2})$ has length 2, 
inducing vanishing of $E_{2}^{p,q}$ for $p\ge 3$. On the other extreme we have
$E_{2}^{p,q}=0$ for $q\ge 2$ since $A$ is hereditary.
\par
Furthermore, the higher derivatives of the invers limit functor can 
be computed as cohomology of the complex

\begin{displaymath}
    \diagram
    End(\F(0))\oplus End(\F(1))
    \oplus End(\F(2))
    \ddto^{
    \matrQQ{\pi^{\ast}}{-\pi_{\ast}}{0}
    {(\pi\circ\epsilon)^{\ast}}{0}{(\pi\circ\epsilon)_{\ast}}
    {0}{\epsilon^{\ast}}{-\epsilon_{\ast}}
    }\\
    \\
    Hom(\F(1),\F(0))
    \oplus Hom(\F(2),\F(0))
    \oplus Hom(\F(2),\F(1))
    \dto^{
    \matrHH{\epsilon^{\ast}}{-1}{\pi_{\ast}}
    }\\
    Hom(\F(2),\F(0))
    \enddiagram
\end{displaymath}
We have $\pi\circ\epsilon=0$ and in fact
$ker \,\pi =im\, \epsilon$. 
It is obvious that 
the second cohomology group vanishes, 
\begin{displaymath}
    \limprojder{2}{Mor(\Lambda_{2})}{Hom_{A}(\F(t),\F(s))}=0
\end{displaymath}
Similar argument proves that  
\begin{displaymath}
    \limprojder{2}{Mor(\Lambda_{2})}{Ext_{A}^1(\F(t),\F(s))}=0
\end{displaymath}
Thus the spectral sequence reduces to the exact sequence
\begin{displaymath}
    0\rightarrow E_{2}^{1,0}\longrightarrow Ext_{\Lambda_{2}}^1(\F,\F)\longrightarrow
    E_{2}^{0,1}\rightarrow 0
\end{displaymath}
and the isomorphism
\begin{displaymath}
    0\rightarrow E_{2}^{1,1}\longrightarrow Ext_{\Lambda_{2}}^2(\F,\F)\rightarrow 0
\end{displaymath}
 
\begin{Prop}\label{Extension tangent space}
    Let $A$ be as above and let
    \begin{displaymath}
\F:\qquad 0\rightarrow W\buildrel{\epsilon}\over{\longrightarrow}
E\buildrel{\pi}\over{\longrightarrow} V\rightarrow 0
    \end{displaymath}
    be a non-split short exact sequence of $A$-modules.
    Suppose $Hom_{A}(V,W)=Hom_{A}(W,V)=0$ and that $End_{A}(V)\simeq
End_{A}(W)\simeq k$
    Then 
    \begin{displaymath}
Ext_{A}^1(E,E)_{0}\simeq ker\{\epsilon^{\ast}\pi_{\ast}:Ext_{A}^1(E,E)\rightarrow 
Ext_{A}^1(W,V)\}
    \end{displaymath}
\end{Prop}

\begin{proof}
The map $Ext_{\Lambda_{2}}^1(\F,\F)\rightarrow Ext_{A}^{1}(E,E))$ factors
through the surjective map $Ext_{\Lambda_{2}}^1(\F,\F)\rightarrow
E_{2}^{0,1}$ and it is enough to consider the image of the map $E_{2}^{0,1}\rightarrow
Ext_{A}^{1}(E,E))$.
    Pick 
    $(a,b,c)\in Ext_{A}^1(W,W)\oplus Ext_{A}^1(E,E)\oplus 
    Ext_{A}^1(V,V)$ such that $\epsilon_{\ast}(a)=\epsilon^{\ast}(b)\in 
    Ext_{A}^1(W,E)$ and $\pi_{\ast}(b)=\pi^{\ast}(c)\in 
    Ext_{A}^1(E,V)$. Thus 
    $\pi_{\ast}\epsilon^{\ast}(b)=\epsilon^{\ast}\pi_{\ast}(b)=0\in 
    Ext_{A}^1(W,V)$ and the image of $E_{2}^{0,1}$ lies in the given kernel.
    \par
    On the other hand, for $b\in ker\{Ext_{A}^1(E,E)\rightarrow 
    Ext_{A}^1(W,V)\}$, we have $\epsilon^{\ast}(b)\in ker\,\pi_{\ast}$ 
    and $\pi_{\ast}(b)\in ker \,\epsilon^{\ast}$. By exactness there exists
$a\in Ext_{A}^1(W,W)$ such 
    that $\epsilon_{\ast}(a)=\epsilon^{\ast}(b)$ and $c\in 
    Ext_{A}^1(V,V)$ such 
    that $\pi^{\ast}(c)=\pi_{\ast}(b)$. This proves that $E_{2}^{0,1}$ maps
surjectively onto the kernel.
\end{proof}
Similarily the obstruction space for deformations under the extension
constrain is given by the image of the map
\begin{displaymath}
Ext_{\Lambda_{2}}^2(\F,\F)\longrightarrow Ext_{A}^2(E,E)
\end{displaymath}
which is trivial since $A$ is hereditary.

Notice that if $Ext_{A}^1(W,V)=0$ and $E$ is a non-trivial extension of $V$ by
$W$, then by the completition theorem \cite{JLS}
there are no simple deformations of $E$, all deformations of $E$ preserve 
the extension structure. This fits with the fact that
$Ext_{A}^1(E,E)_{0}\simeq Ext_{A}^1(E,E)$ and there are no obstructions.
The map $\epsilon^{\ast}\pi_{\ast}$ of Proposition \ref{Extension tangent
space} is a composition of two surjective maps, and therefore itself
surjective. It follows that the moduli space of extensions is smooth of 
dimension 
\begin{displaymath}
    dim\,Ext_{A}^1(E,E)-dim\,Ext_{A}^1(W,V)
\end{displaymath}

\begin{Prop}
    Let $\F$ be an non-split exact sequence as above, with dimension vector $\beta$
and 
    $\gamma$ for $W$ and $V$ respectively. Then
    \begin{displaymath}
dim\,Ext_{\Lambda_{2}}^1(\F,\F)=1-\langle 
\beta+\gamma,\beta+\gamma\rangle+\langle\beta,\gamma\rangle
    \end{displaymath}
\end{Prop}
\begin{proof}
    Follows immedeately from the dimension formula and 
    \begin{displaymath}
    dim\,Ext^1_{\Lambda_{2}}(\F,\F)=dim\,Ext^1_{A}(E,E)_{0}
    =dim\,Ext_{A}^1(E,E)-dim\,Ext_{A}^1(W,V)
    \end{displaymath}
\end{proof}

Notice that the
codimension of the non-simple locus at $\F$ is given by $-\langle\beta,\gamma\rangle$
where $\beta+\gamma=\alpha$.

\subsection{An example}\label{An example}
Consider the dimension vector $\alpha=(2,2,2;3,3)$, $n=6$. Then
$dim\,Simp_{\alpha}(A)=1+6^2-(2^2+2^2+2^2)-(3^2+3^2)=7$. The minimum value of
the codimension $-\langle\beta,\gamma\rangle$ is obtained for
$\beta=(2,2,1;3,2)$, $\gamma=(0,0,1;0,1)$ which gives
$-\langle\beta,\gamma\rangle=2$. This shows that $\M_{\alpha}$ in this case is a
7-dimensional smooth variety with a 5-dimensional closed subvariety of non-simple
modules.

\begin{Rmk}
Notice that the non-simple locus does not have to be a smooth subvariety, even if
the obstruction space of deformations under the extension constrain is trivial. A
non-simple module can be written as an extension in many different ways.
\end{Rmk}

\section{Maximally iterated extensions}\label{Maximally iterated extension}
\subsection{Maximally iterated extensions}\label{MIE}

Let $M$ be a left $A$-module. The structure map $\rho:A\rightarrow End_{k}(M)$
of the module $M$ is determined by two $n\times n$-matrices $X=\rho(x)$ and $Y=\rho(y)$,
up to
simultanous conjugation, satisfying $X^3=Y^2=I$. The eigenspaces of $X$ and $Y$
corresponding to the three cube roots and the two square roots of unity,
denoted $X_{i}$, $i=1,2,3$, resp. $Y_{j}$, $j=1,2$, have dimensions given by the
dimension vector $\alpha=(\alpha_{1},\alpha_{2},\alpha_{3};\alpha^{\prime}_{1},
\alpha^{\prime}_{2})$.  
\par
Let 
\begin{displaymath}
M=M_{r}\supset M_{r-1}\supset \dots\supset M_{1}\supset M_{0}=0
\end{displaymath}
be any composition series of $M$. We say that $M$ is a maximally iterated
extension if 
all factors $\overline{M}_{i}=M_{i}/M_{i-1}$ have dimension 1.  

\subsection{Counting Maximally iterated extensions}\label{Couting MIE}

There exists a set of distinct points in 
$\M_{\alpha}$, corresponding to the weak classes of maximally iterated extensions
of  $A$-modules. 

\begin{Prop}
Let $\alpha$ be a Schur root. 
The set of maximally iterated extensions of dimension vector $\alpha$ is 
    finite of cardinality $N=\frac{1}{2}(d+n+1)$.
\end{Prop}

\begin{proof}
The associated semi-simple module $\overline{M}=\oplus_{i=1}^n\overline{M}_{i}$
of a maximally iterated extension $M$ can be written as a sum
\begin{align*}
\overline{M}=k(1,1)^{a_{1}}&\oplus k(\omega,1)^{a_{2}}
\oplus k(\omega^2,1)^{a_{3}}\\
&\oplus k(1,-1)^{\alpha_{1}-a_{1}}
\oplus k(\omega,-1)^{\alpha_{2}-a_{2}}\oplus k(\omega^2,-1)^{\alpha_{3}-a_{3}}
\end{align*}
for non-negative integers $a_{1}$, $a_{2}$, $a_{3}$, satisfying
$a_{i}\le \alpha_{i}$, $i=1,2,3$. The 
    generating function for the problem of finding all integers fullfilling
this condition is given by
    \begin{align*}
g(x)&=(1+x+x^2+\ldots+x^{\alpha_{1}})(1+x+x^2+\ldots+x^{\alpha_{2}})
(1+x+x^2+\ldots+x^{\alpha_{3}})\\
&=\frac{1-x^{\alpha_{1}+1}}{1-x}\frac{1-x^{\alpha_{2}+1}}{1-x}
\frac{1-x^{\alpha_{3}+1}}{1-x}\\
&=(1-x^{\alpha_{1}+1})(1-x^{\alpha_{2}+1})(1-x^{\alpha_{3}+1})
\frac{1}{(1-x)^3}\\
&=(1-x^{\alpha_{1}+1})(1-x^{\alpha_{2}+1})(1-x^{\alpha_{3}+1})
\sum_{i=0}^{\infty}
\vect{i+2}2x^i\\
    \end{align*}
    Let $m=a_{1}+a_{2}+a_{3}$. The coefficient of the $x^m$-term is given by 
    \begin{align*}
N=&\vect{m+2}2-\vect{m-\alpha_{1}+1}2-\vect{m-\alpha_{2}+1}2-
\vect{m-\alpha_{3}+1}2\\
&+\vect{m-\alpha_{1}-n_{2}}2+\vect{m-\alpha_{1}-\alpha_{3}}2+
\vect{m-\alpha_{2}-\alpha_{3}}2\\
&-\vect{m-\alpha_{1}-\alpha_{2}-\alpha_{3}}2\\
    \end{align*}
    By definition $m=\alpha_{1}^{\prime}$ and $n-m=\alpha_{2}^{\prime}$. The
conditions for existence of a Schur root give
    \begin{displaymath}
\alpha_{i}+m=\alpha_{i}+\alpha_{1}^{\prime}\le \alpha_{i}+\alpha_{j}+\alpha_{k}\quad 
\Rightarrow \quad m\le \alpha_{j}+\alpha_{k}
    \end{displaymath}
    and
    \begin{displaymath}
\alpha_{i}+\alpha_{1}+\alpha_{2}+\alpha_{3}-m
=\alpha_{i}+\alpha_{2}^{\prime}\le 
\alpha_{i}+\alpha_{j}+\alpha_{k}\quad
\Rightarrow
\quad 
m\ge \alpha_{i}
    \end{displaymath}
    for $\{i,j,k\}=\{1,2,3\}$.
    \newline
    Thus $m-\alpha_{i}-\alpha_{j}\le 0$ and $m-\alpha_{i}+1\ge 1$ and we get
    \begin{align*}
N&=\vect{m+2}2-\vect{m-\alpha_{1}+1}2-\vect{m-\alpha_{2}+1}2-
\vect{m-\alpha_{3}+1}2\\
&=\frac{1}{2}\big((m^2+3m+2)-(m^2-2m\alpha_{1}+m+{\alpha_{1}}^2-\alpha_{1})\\
&\quad -(m^2-2m\alpha_{2}+m+{\alpha_{2}}^2-\alpha_{2})-
(m^2-2m\alpha_{3}+m+{\alpha_{3}}^2-\alpha_{3})\big)\\
&=\frac{1}{2}(-2m^2+2+2mn-\sum_{i} {\alpha_{i}}^2+n)\\
&=m(n-m)+\frac{1}{2}(n-\sum _{i}{\alpha_{i}}^2)+1\\
    \end{align*}
    Changing back to $m=\alpha^{\prime}_{1}$, $n-m=\alpha^{\prime}_{2}$, this
takes the form
    \begin{displaymath}
N=\alpha^{\prime}_{1}\alpha^{\prime}_{2}+\frac{1}{2}(n-\sum _{i}{\alpha_{i}}^2)+1
    \end{displaymath}
    But $n=\alpha^{\prime}_{1}+\alpha^{\prime}_{2}$, so 
    $\alpha^{\prime}_{1}\alpha^{\prime}_{2}=\frac{1}{2}(n^2-\sum_{j}
(\alpha^{\prime}_{j})^2)$ 
    and we get
    \begin{displaymath}
N=\frac{1}{2}(n+n^2)-\frac{1}{2}\sum _{i}{\alpha_{i}}^2-\frac{1}{2}\sum_{j}
(\alpha^{\prime}_{j})^2+1
    \end{displaymath}
    From Theorem \ref{Westbury} we know that 
    \begin{displaymath}
d=dim\,Ext_{A}^1(E,E)=n^2-\sum_{i} {\alpha_{i}}^2-\sum_{j} (\alpha^{\prime}_{j})^2+1
    \end{displaymath}
    showing that 
    \begin{displaymath}
N=\frac{1}{2}(n+n^2)-\frac{1}{2}n^2-\frac{1}{2}+\frac{d}{2}+1
=\frac{1}{2}(n+d+1)
    \end{displaymath}
\end{proof}
Notice that a consequence of this formula is that the number $n+d$ must be odd
for any dimension vector $\alpha$.

\section{Parametrisation of Maximally iterated
extensions}\label{Parametrisation}

\subsection{Strong vs. weak equivalence relation}\label{The dimension theorem}

 A point in the subvariety
$\D_{\alpha}$ corresponds to an equivalence class of
indecomposable representations, all with conjugate composition factors.
Under the stronger conjugation equivalence relation this class splits up in
several conjugacy classes. Since indecomposable modules deform into simples,
but not vice versa, the most extreme modules are the maximally iterated
extensions, where all composition factors are of dimension 1.
Let 
$\overline{E}\in \D_{\alpha}$
be the weak equivalence class of a maximally iterated extension $E$. The problem
is to understand how this equivalence class splits into different
classes under the stronger conjugation equivalence relation.
\par
More precisely, given a dimension vector $\alpha$ we choose one of the
$\frac{1}{2}(n+d+1)$ points of $\M_{\alpha}$ of the last section,
corresponding to a maximally iterated extension $E$. The corresponding
semi-simple module $\overline{E}$ defines a set of 1-dimensional simple
representations, possible of multiplicity greater than 1. Let $\Gamma$ be a
representation graph for this set, i.e. an ordering of the given simple
factors. Following Laudal (\cite{L3}) we define $Ind_{\Gamma}(A)$ as the
non-commutative scheme of indecomposable $\Gamma$-representations. Our task is
to describe this non-commutative scheme for the modular group ring.

\subsection{Some usefull lemmas}\label{Some usefull lemmas}
For an upper triangular $n\times n$-matrix $Z=(z_{ij})$ we introduce the notation
\begin{displaymath}
i\mathop{\sim}_{Z}j
\end{displaymath}
to indicate that the diagonal elements $z_{ii}=z_{jj}$. This is obviously an
equivalence relation on the set $\{1,\dots,n\}$. The subscript $Z$, referring to
the matrix, will be omitted if the reference to the matrix is unambiguous.

\begin{Lem}\label{Jordan}
Let $X$ be an upper triangular $n\times n$-matrix on Jordan block form,
satisfying the equation $X^3=I$. Then $X$ is diagonal.
\end{Lem}

\begin{proof}
Suppose $X=(x_{ij})$ is on Jordan form, i.e. the three eigenvalues on the
diagonal is sequentially ordered and all off-diagonal
entries are zero, except possible for $x_{i,i+1}$ if $x_{ii}=x_{i+1,i+1}$. But
$(X^3)_{i,i+1}=3(x_{ii})^2x_{i,i+1}=0$, showing that $X$ is diagonal.
\end{proof}

\begin{Lem}\label{utm}
Let $Z$ be a diagonalisable upper triangular $n\times n$-matrix. Then there
exists an upper triangular invertible matrix 
$U$ such that $U^{-1}ZU$ is diagonal.
\end{Lem}

\begin{proof}
Let $P$ be an invertible $n\times n$-matrix. By the LU-decomposition
Theorem $P=\Gamma LU$ where $\Gamma$ is a permutation matrix, $L$ is lower
triangular and $U$ is upper triangular.
Suppose $P^{-1}$ diagonalize $Z$, i.e. $PZP^{-1}=D$, $D$ is diagonal. Then 
\begin{displaymath}
UZU^{-1}=L^{-1}\Gamma^{-1}D\Gamma L
\end{displaymath}
The left-hand side of this equation is an upper triangular matrix, while the 
right-hand side is lower triangular since conjugation of a diagonal matrix by
a permutation matrix again is diagonal.
Thus the matrix $UZU^{-1}$ is diagonal.
\end{proof}

The outcome of these two lemmas is that given a pair of upper triangular
matrices $(X,Y)$, where $X$ satisfies the equation $X^3=I$, we can always
diagonalise $X$, without disturbing the upper triangular form of $Y$. 
\par
The next lemma gives an explicit form for  the upper triangular matrices $Y$,
satisfying $Y^2=I$. 
\par
Let $\P(i,k)$, $i<k$ be given by 
\begin{displaymath}
\{i=i_{0}<i_{1}<\dots<i_{m}<i_{m+1}=k\}
\in \mathcal{P}(i,k)
\end{displaymath}
if $m$ is odd and
$y_{i_{r}i_{r}}=(-1)^ry_{ii}$ for $r=0,1,2,\dots,m+1$. Define the sequence
$\chi:\mathrm{Z}_{+}\rightarrow\mathrm{Z}_{+}$ recursively by $\chi(1)=1$ and
$\chi(n)=\sum_{i=1}^{n-1}\chi(i)\chi(n-i)$.
\begin{Lem}\label{2rot}
Let $Y=(y_{ij})$ be an upper triangular $n\times n$-matrix, satisfying
$Y^2=I$. Then $y_{ii}=\pm 1$ and for $y_{ii}= y_{kk}$, $i<k$ we have 
\begin{displaymath}
y_{ik}=\sum_{\mathcal{P}(i,k)}(-y_{ii})\frac{\chi(\nu)}{2^{m}}\prod_{s=0}^{m}
y_{i_{s}i_{s+1}}
\end{displaymath}
where $m=2\nu +1$.
The entry $y_{ik}$ can be
choosen freely  whenever $y_{ii}\neq y_{kk}$. 
\end{Lem}

An immedeate consequence of this lemma is that the space of upper
triangular $n\times n$-matrices, satisfying $Y^2=I$ with eigenspaces of
dimension $\alpha^{\prime}_{1},\alpha^{\prime}_{2}$ is an affine space
of dimension $\alpha_{1}^{\prime}\cdot \alpha_{2}^{\prime}=\frac{1}{2}
(n^2-({\alpha_{1}}^{\prime})^2-({\alpha_{2}}^{\prime})^2)$.

\begin{proof}
Since $Y$ is upper triangular and $Y^2=I$ the diagonal entries obviously
satisfy the equation $y_{ii}^2=1$.
\par
The relation $Y^2=I$ can be written as a system of $\vect{n}2$ quadratic equations
in the $\vect{n}2$ off-diagonal entries $y_{ik}$, $1\le i<k\le n$
\begin{equation}
(y_{ii}+y_{kk})y_{ik}+\sum_{j=i+1}^{k-1} y_{ij}y_{jk}=0
\end{equation}
The proof of the lemma has two parts. First we show that the given vector $(y_{ik})$
is a
solution of the system (1). Then we show that all solutions are of this
form. We proceed by induction on the difference $\delta=k-i$. If $\delta=1$ 
the defining equation reduces to 
\begin{displaymath}
(y_{ii}+y_{i+1,i+1})y_{i,i+1}=0
\end{displaymath}
If $y_{ii}+y_{i+1,i+1}=0$ there are no conditions on $y_{i,i+1}$. If
$y_{ii}+y_{i+1,i+1}\neq 0$ it follows that $y_{i,i+1}=0$. This fits with the
given formula since in that case $\P(i,i+1)=\emptyset$. Thus the lemma is true
for $\delta=1$.
\par
Now suppose the formula is valid for $\delta<t$ and let $(i,k)$ satisfy $k-i=t$.
If $y_{ii}+y_{kk}=0$ we have to show that 
\begin{displaymath}
\sum_{j=i+1}^{k-1} y_{ij}y_{jk}=0
\end{displaymath}
Since $y_{ii}\neq y_{kk}$ the sum splits into two parts
\begin{align*}
\sum_{j=i+1}^{k-1} y_{ij}y_{jk}
&=\sum_{i\mathop{\sim}j} y_{ij}y_{jk}+\sum_{j\mathop{\sim}k} y_{ij}y_{jk}\\
&=\sum_{i\mathop{\sim}j}( \sum_{\mathcal{P}(i,j)}(-y_{ii})\frac{\chi(p)}{2^{p}}
\prod_{s=0}^{p}y_{i_{s}i_{s+1}})y_{jk}\\
&\hspace{20mm}+\sum_{j\mathop{\sim}k}
y_{ij}(\sum_{\mathcal{P}(j,k)}(-y_{jj})\frac{\chi(q)}{2^{q}}
\prod_{s=0}^{q}y_{j_{t}j_{t+1}})
\end{align*}
For each sequence 
\begin{displaymath}
i=i_{0}<i_{1}<... < i_{m+1}=k
\end{displaymath}
there are two apperances of the term $y_{i_{0}i_{1}}y_{i_{1}i_{2}}\dots
y_{i_{m-1}k}$, one in each sum. In the first
sum $j=i_{m}$ and in the second sum $j=i_{1}$. Thus the two terms add up to
\begin{align*}
(-y_{ii})\frac{\chi(m-1)}{2^{m-1}}
\prod_{s=0}^{m-1}& y_{i_{s}i_{s+1}}y_{jk}+y_{ij}(-y_{jj})
\frac{\chi(m-1)}{2^{m-1}}
\prod_{s=1}^{m}y_{j_{t}j_{t+1}}\\
&=
 (-y_{i_{0}i_{0}})\frac{\chi(m-1)}{2^{m-1}}
 \prod_{s=0}^{m-1} y_{i_{s}i_{s+1}}y_{i_{m}k}
 +y_{ii_{1}}(-y_{i_{1}i_{1}})\frac{\chi(m-1)}{2^{m-1}}
 \prod_{t=1}^{m}y_{i_{t}i_{t+1}}\\
& =
 (-y_{i_{0}i_{0}}-y_{i_{1}i_{1}})\frac{\chi(m-1)}{2^{m}}
 \prod_{t=0}^{m-1}y_{i_{t}i_{t+1}}=0
\end{align*}
\par
For $y_{ii}=y_{kk}=\pm 1$ the equation (1) takes the form
\begin{displaymath}
y_{ik}=-\frac{1}{2}y_{ii}\sum_{j=i+1}^{k-1} y_{ij}y_{jk}
\end{displaymath}
and we have to show that 
\begin{displaymath}
-\frac{1}{2}y_{ii}\sum_{j=i+1}^{k-1} y_{ij}y_{jk}=
\sum_{\mathcal{P}(i,k)}(-y_{ii})\frac{\chi(m-1)}{2^{m-1}}\prod_{s=0}^{m}
y_{i_{s}i_{s+1}}
\end{displaymath}
We need a lemma.

\begin{Lem}\label{P}
Let $i<j<k$ and suppose $y_{ii}=y_{kk}$. Then
\begin{displaymath}
[\bigcup_{i\mathop{\sim}j}(\P(i,j)\times \P(j,k))]\cup \{i<j<k\,\vert\,
y_{ii}\neq y_{jj}\}= \P(i,k)
\end{displaymath}
\end{Lem}

\begin{proof}
Let $\I=\{i=i_{0}<i_{1}<\dots <i_{p}=j\}\in \P(i,j)$ and 
$\J=\{j=j_{0}<j_{1}<\dots <j_{q}=k\}\in \P(j,k)$. Define the join
\begin{align*}
\I\cdot\J&=\{i=i_{0}<i_{1}<\dots <i_{p}=j_{0}<j_{1}<\dots j_{q}=k\}\\
&=\{i=i_{0}<i_{1}<\dots <i_{p}<i_{p+1}<\dots i_{p+q}=k\}
\end{align*}
where we have put $i_{p+t}=j_{t}$ for $t=0,1,\dots , q$. 
 If $p,q$ are even, so is $p+q$, and for $r>p$ we have 
\begin{displaymath}
y_{i_{r}i_{r}}=y_{j_{r-p}j_{r-p}}=(-1)^{r-p}y_{jj}
=(-1)^{r-p}y_{i_{p}i_{p}}=(-1)^{r-p+p}y_{i_{0}i_{0}}
=(-1)^ry_{ii}
\end{displaymath}
Thus $\I\cdot\J\in\P(i,k)$. 
The set $\{i<j<k\,\vert\,y_{ii}\neq y_{jj}\}$ is easily seen to be in
$\P(i,k)$.
\par
Now pick some $\K\in\P(i,k)$ and choose $j\in \K$ such that
$y_{ii}=y_{jj}=y_{kk}$. Then $\K=\I\cdot\J$ where
\begin{displaymath}
\I=\{i=k_{0}<\dots<k_{p}=j\}\quad\textrm{and}\quad
\J=\{j=k_{p}<k_{p+1}\dots<k_{p+q}=k\}
\end{displaymath}
and $y_{ii}=y_{k_{p},k_{p}}=(-1)^{p}y_{ii}$ is possible only if $p=j-i$ is an
even number.
\end{proof}

\textit{Proof of Lemma \ref{2rot} continued.}\\
Our assumption is that 
\begin{displaymath}
y_{ik}=\sum_{\mathcal{P}(i,k)}(-y_{ii})\frac{\chi(\nu)}{2^{m}}\prod_{s=0}^{m}
y_{i_{s}i_{s+1}}
\end{displaymath}
for $k-i=\delta<t$. For $k-i=t$ we get 
\begin{align*}
(y_{ii}&+y_{kk})y_{i k}+\sum_{j=i +1}^{k-1}y_{ij}y_{jk}\\
&=(y_{ii}+y_{kk})y_{i k}+\sum_{i\mathop{\not\sim}j}y_{ij}y_{jk}\\
&+
\sum_{j=i +1}^{k-1}
(\sum_{\mathcal{P}(i,j)}(-y_{ii})\frac{\chi(n-1)}{2^{n-1}}\prod_{s=0}^{n}
y_{i_{s}i_{s+1}})
(\sum_{\mathcal{P}(j,k)}(-y_{jj})\frac{\chi(m-1)}{2^{m-1}}\prod_{t=0}^{m}
y_{i_{t}i_{t+1}})\\
&=(y_{ii}+y_{kk})y_{i k}+\sum_{i\mathop{\not\sim}j}y_{ij}y_{jk}\\
&+\sum_{j=i +1}^{k-1}
(\sum_{\P(i,j)\times\P(j,k)}(-y_{ii})\frac{\chi(n-1)\chi(m-1)}
{2^{n+m-2}}\prod_{s=0}^{n+m}
y_{i_{s}i_{s+1}})\\
&=(y_{ii}+y_{kk})y_{i k}+\sum_{i\mathop{\not\sim}j}y_{ij}y_{jk}
+\sum_{j=i +1}^{k-1}
(\sum_{\P(i,k)}(-y_{ii})\frac{\chi(n+m-2)}
{2^{n+m-2}}\prod_{s=0}^{n+m}
y_{i_{s}i_{s+1}})
\end{align*}
where $s_{n+j}=t_{j}$ for $j=0,1,\dots, m$.
\par
For $y_{ii}+y_{kk}\neq 0$ we get
\begin{displaymath}
y_{i k}=-\sum_{j=i +1}^{k-1}
(\sum_{\P(i,k)}(-y_{ii})\frac{\chi(u-1)}{2^{u-1}}\prod_{s=0}^{u}
y_{i_{s}i_{s+1}})
\end{displaymath}
and for $y_{ii}+y_{kk}= 0$ this relation gives no new condition on $y_{ik}$.
\end{proof}

In the next lemma we establish the explicit form for the stabilising
$GL_{n}$-subgroup of the diagonal matrix $X$, satisfying $X^3=I$, with
eigenspaces of dimension $\alpha_{1},\alpha_{2},\alpha_{3}$.
 
\begin{Lem}\label{dimstab}
The stabiliser subgroup $G_{X}\subset GL_{n}$, stabilising the diagonal matrix $X$,
is given by
\begin{displaymath}
G_{X}=\{(g_{ij})\in GL_{n}\,\vert\, g_{ij}(x_{ii}-x_{jj})=0\}
\end{displaymath}
It is a linear group of dimension
\begin{displaymath}
dim\,G_{X}=\alpha_{1}^2+\alpha_{2}^2+\alpha_{3}^2
\end{displaymath}
isomorphic to the product group $M_{\alpha_{1}}(k)\times
M_{\alpha_{2}}(k)\times M_{\alpha_{3}}(k)$.
\end{Lem}

\begin{proof}
Let $g\in G_{X}$. The relation $gX=Xg$ gives for the entries 
$g_{ij}x_{jj}=x_{ii}g_{ij}$ or
$g_{ij}(x_{jj}-x_{ii})=0$. Thus we must have $g_{ij}=0$ for all pairs $(i,j)$
such that $x_{jj}-x_{ii}\neq 0$, i.e. $g_{ij}=0$ whenever $x_{ii}\neq x_{jj}$.
There are
$n^2$ entries in the whole matrix, $n$ are on the diagonal and 
$\sum_{i=1}^3{\alpha_{i}}^2-\alpha_{i}$ entries correspond to off-diagonal pairs
$(i,j)$ with $x_{ii}=x_{jj}$, adding up to the given dimension.
\par
Since $X$ is diagonal with at most three different eigenvalues, conjugation by
a certain permutation matrix shows that the stabiliser subgroup is isomorphic
to the product of three full matrix groups, of given dimension.

\end{proof}

Denote by 
\begin{displaymath}
\Y(\Gamma)=\{Y\in M_{n}(k)\,\vert
\,Y\,\textrm{upper}\,\textrm{triangular},\,Y^2=I, \,Y\in\Gamma_{Y}\}
\end{displaymath}
By $Y\in\Gamma_{Y}$ we mean that $y_{11}, ... ,y_{nn}$ is choosen in
accordance with the representation graph $\Gamma$.

We have $dim\,\Y(\Gamma)=\alpha^{\prime}_{1}\alpha_{2}^{\prime}
=\frac{1}{2}(n^2-\sum_{j}(\alpha_{j}^{\prime})^2)$
with a natural basis given by $\{(i,j)\,\vert\,i<j,\,i\mathop{\not\sim}_{Y}j\}$. 
The subgroup of $G_{X}$ for $X\in\Gamma_{X}$ stabilising $\Y(\Gamma)$ is denoted
$G_{X,\Y(\Gamma)}$. 

Recall that the original problem of this section was the following: Given a
weak equivalence class
$\overline{E}\in \D_{\alpha}$
of a maximally iterated extension $E$. This equivalence class splits up in a
bunch of equivalence classes under the stronger conjugation relation,
described by the non-commutative scheme $Ind_{\Gamma}(A)$.
Then we have the following

\begin{Prop}
 Then there is a one-to-one correspondance 
\begin{displaymath}
Ind_{\Gamma}(A))\longleftrightarrow \Y(\Gamma)/G_{X,\Y(\Gamma)}
\end{displaymath}
\end{Prop}

\begin{proof}
A representation in $Ind_{\Gamma}(A)$ is given by two upper triangular matrices
$X$ and $Y$, satisfying $X^3=Y^2=I$. By Lemma \ref{Jordan} we can assume that
$X$ is diagonalisable and by Lemma \ref{utm} we can in fact assume that $X$ is
diagonal and $Y$ is upper triangular. The parameter space $\Y(\Gamma)$ is
described in Lemma \ref{2rot}. 
\end{proof}

\section{Relation to modular forms}\label{Modular Forms}

The algebra of modular forms is freely generated as a commutative ring by the
the two forms $E_{4}$ and $E_{6}$. Let $f(n)$ be the dimension of
the space of
modular forms of dimension $2n$, $n\ge 2$. The generating function for $f(n)$
is given by
\begin{displaymath}
\sum_{n=0}^{\infty} f(n)x^{n}=\frac{1}{1-x^2}\cdot\frac{1}{1-x^3}
\end{displaymath}
\par
We have seen that $dim\,Simp_{n}(A)=dim\, Simp_{\alpha}(A)$ for
$\vert\alpha\vert=n$, where the partitions
$n=\alpha_{1}+\alpha_{2}+\alpha_{3}=\alpha_{1}^{\prime}+\alpha_{2}^{\prime}$
are at most evenly distributed. Let $\alpha=\beta+\gamma$. We are interested in
the number
$-\langle\beta,\gamma\rangle
=\vert\beta\vert\cdot\vert\gamma\vert-\sum_{i=1}^3\beta_{i}\gamma_{i}
-\sum_{j=1}^2\beta_{j}^{\prime}\gamma_{j}^{\prime}$. It can be proven that the
minimal value of this number is obtained when 
$\vert\gamma\vert=1$ (or $\vert\beta\vert=1$).

For fixed
$\vert\beta\vert$ and $\vert\gamma\vert$, the smallest value of
$-\langle\beta,\gamma\rangle$ is again achieved
when $\beta$ and $\gamma$ are at most evenly distributed. On the other hand,
$\vert\beta\vert\cdot\vert\gamma\vert$ has its minimal positive value when one of
the
factors equals 1. Thus the minimal value for $-\langle\beta,\gamma\rangle$ is
given by $n-1-\beta_{i}-\beta_{j}^{\prime}$ for maximal
$\beta_{i}=\alpha_{i}-1$ and $\beta_{j}^{\prime}=\alpha_{j}^{\prime}-1$. Since
$\alpha$ is at most evenly distributed we have 
$\beta_{i}=[\frac{n-1}{3}]$ and  $\beta_{j}^{\prime}=[\frac{n-1}{2}]$, where
$[q]$ denotes the greatest integer, less than $q$.

 The generating function for the number  $[\frac {n-1}{k}]$ is given by
\begin{displaymath}
(x+x^2+x^3+\dots)(x^k+x^{2k}+x^{3k}+\dots)=\frac{x^{k+1}}{(1-x)(1-x^k)}
\end{displaymath}
Thus the generating function for the sequence
$n-1-[\frac{n-1}{3}]-[\frac{n-1}{2}]$ is given by

\begin{align*}
P(x)&=x(x+2x^2+3x^3+\dots)-\frac{x^4}{(1-x)(1-x^3)}-\frac{x^3}{(1-x)(1-x^2)}\\
&=\frac{x^2}{(1-x)^2}-\frac{x^4}{(1-x)(1-x^3)}-\frac{x^3}{(1-x)(1-x^2)}\\
&=\frac{1}{(1-x^2)(1-x^3)}-1
\end{align*}
Thus except for the one of weight 0, the number of modular
forms of weight $2n$ equals the codimension of the non-simple stratum of the
highest-dimensional component of the variety of semi-simple representations of
$A$ of dimension $n$.

\end{document}